\documentclass[11pt]{article}
\usepackage{amsfonts,amsmath,amssymb,amsthm, amscd}
\usepackage[dvips]{epsfig,graphics}

\numberwithin{equation}{section}
\newtheorem{theorem}{Theorem}[section]
\newtheorem{prop}[theorem]{Proposition}

\newtheorem{conj}[theorem]{Conjecture}

\theoremstyle{definition}
\newtheorem{definition}[theorem]{Definition}

\newtheorem{remark}[theorem]{Remark}

\newcommand{\D}{\Delta}

\def\<{{\langle}}
\def\>{{\rangle}}

\def\S{{\Sigma}}

\def\sx{\sigma_x}

\begin{document}

\title{Twisted Alexander Polynomials and Representation Shifts}

\author{Daniel S. Silver \and Susan G. Williams\thanks{Both authors partially supported by NSF grant
DMS-0304971.} \\ {\em
{\small Department of Mathematics and Statistics, University of South Alabama}}}

%date August 2, 2007
\maketitle %{\setlength{\linewidth}{2in}

\begin{abstract}\noindent For any knot, the following are equivalent. (1) The infinite cyclic cover has uncountably many finite covers; (2) there exists a finite-image representation of the knot group for which the twisted Alexander polynomial
vanishes; (3) the knot group admits a finite-image representation such that the image of the fundamental group of an incompressible Seifert surface is a proper subgroup of the image of the commutator subgroup of the knot group. \end{abstract}

\noindent {\it Keywords:} Knot, knot group, twisted Alexander polynomial, representation shift
\begin{footnote}{Mathematics Subject Classification:  
Primary 57M25; secondary 37B40.}\end {footnote}

%%%%%%%%%%%%%%%%%%%%%%%%%%%%%% 1. INTRODUCTION %%%%%%%%%%%%%
\section{Introduction} Fibered knots have been studied in settings of complex algebraic geometry and dynamical systems as well as topology. Their relative simplicity accounts in part for their appeal. Among the ``simplest hyperbolic knots" catalogued by  P.  Callahan, J. Dean and J. Weeks \cite{cdw}, those hyperbolic knots with 6 or fewer ideal tetrahedra in their complements, an overwhelming majority are fibered. 

Each of the three conditions listed in the abstract is known to imply that the knot $k$ is not fibered.  We conjecture that conversely, each of these conditions holds for every nonfibered knot.  That (3) implies (2) we derive from a result of Friedl and Vidussi \cite{fvJAMS}.  The equivalence of (1) and (3) follows from our earlier work \cite{swJPAA}.
In \cite{swJPAA} we also stated a variant of the above conjecture, given below as Conjecture \ref{swconj}.  The weaker Conjecture \ref{folkconj}, which is unattributed, is closely related to a conjecture in \cite{fv2}.
%J. Hempel proved that knot groups are residually finite \cite{hempel}, meaning that any pair of distinct elements can be distinguished in  some finite quotient group. In view of this fact, the authors conjectured a knot fibering criterion based on finite-image representations of the commutator subgroup of the knot group (Conjecture \ref{swconj} below). Another conjecture, of uncertain origin, is based on twisted Alexander polynomials associated to finite-image representations of the knot group.
%(see Conjecture \ref{folkconj}). 

%

% A weak form of subgroup separability would imply both conjectures (See section 3). 

Section 2 presents basic notions of representation shifts and twisted Alexander polynomials. The main result appears in Section 3. We relate Conjectures \ref{swconj} and \ref{folkconj} in the conclusion. 

Many of our results about knot complements generalize in a straightforward manner for arbitrary $3$-manifolds. 

We are most grateful to Stefan Friedl for insightful comments and suggestions. 

\section{Fibered knots, rep shifts and twisted homology}
We recall that a knot $k\subset {\mathbb S}^3$ is fibered if its exterior admits a locally trivial fibration over the circle. By combined theorems of Neuwirth \cite{neuwirth} and Stallings \cite{stallings}, this condition is equivalent to the requirement that the commutator subgroup $\pi'$ of the knot group $\pi = \pi_1({\mathbb S}^3 \setminus k)$ is finitely generated. Then $\pi'$ is in fact isomorphic to the fundamental group of a minimal-genus Seifert surface $S$ for $k$.

If $k$ is fibered, then $\pi'$, being finitely generated, has countably many subgroups of finite index. Equivalently, the infinite cyclic cover $X_\infty$ of $k$ has countably many finite covers. However, for nonfibered knots the situation is different. For every known example, the number of finite covers is uncountable. The authors conjectured this for all nonfibered knots \cite{swTAMS}, and proved it for nonfibered knots of genus 1 \cite{sw07}.

For any knot, the finite-index subgroups of $\pi'$ can be found using representation shifts. We sketch the basic idea, referring the reader to 
\cite{swIJM} \cite{swJPAA} \cite{swTAMS} for details. 

Given a finite group $\S$, consider the space $\Phi_\S= {\rm Hom}(\pi', \S)$ of homomorphisms $\rho: \pi' \to \S$. The topology is the compact-open topology, where $\pi'$ and $\S$ are discrete. Conjugation in $\pi$ by a meridian $x$ induces a homeomorphism $\sx$ described by 
$$\sx \rho(a) = \rho(x^{-1}a x),\ \forall a \in \pi'.$$ 
The pair $(\Phi_\S, \sx)$, called a {\it representation shift}, is a dynamical system that, up to topological conjugacy, is independent of the choice of meridian. Consequently, the topological entropy $h_\S$, a measure of complexity of $\sx$, is an invariant of $k$ and $\S$. 

We say $\rho\in\Phi_\S$ has {\it period} $r$ if $\sx^r\rho=\rho$.  The set ${\rm Fix}\ \sx^r$ of points with period $r$ can be identified with ${\rm Hom}(\pi_1 M_r, \S)$, where $M_r$ denotes the $r$-fold cyclic cover of ${\mathbb S}^3$ branched over $k$; in fact, when $\S$ is abelian, 
${\rm Fix}\ \sx^r$ and ${\rm Hom}(\pi_1 M_r, \S)$ are isomorphic as abelian groups. 
The following proposition is proved in \cite{swTAMS} (see also \cite{sw07} for a sketch of the proof).  

\begin{prop} \label{summary} For any knot $k$,  the following  statements are equivalent.
\item{(1)} The infinite cyclic cover of $k$ has  uncountably many finite covers. 
\item{(2)} The representation shift $\Phi_\S$ is uncountable, for some finite group $\S$. 
\item{(3)} The topological entropy $h_{\S}$ is positive, for some finite group $\S$.
\item{(4)} $\displaystyle \lim_{r\to \infty} \displaystyle{1\over r} \log |{\rm Hom}(\pi_1 M_r, \S)|$ is positive, for some finite group $\S$. 
\end{prop}

The following conjecture of \cite{swTAMS} proposes a characterization of fibered knots. 
\begin{conj}\label{swconj}

A knot $k$ is fibered iff $h_\S =0$ for every finite group $\S$.

\end{conj}

\begin{remark} As noted above, the forward implication is known. Conjecture \ref{swconj} can be restated as: $k$ is nonfibered iff any of the equivalent conditions of Proposition \ref{summary} holds. \end{remark}

It is well known that the Alexander polynomial $\D_k$ of a fibered knot is a monic polynomial with degree equal to twice the genus of the knot. While useful, such a fibering obstruction is 
far from complete. The Alexander polynomial of any fibered knot is also the Alexander polynomial of infinitely many nonfibered knots.
Twisted Alexander polynomials, introduced by X.S. Lin \cite{lin} in 1990, provide more sensitive fibering obstructions. 

The Alexander polynomial $\D_k$ is the $0$th characteristic polynomial of  $H_1 X_\infty$, regarded as a ${\mathbb Z}[t, t^{-1}]$-module. The module has a square presentation matrix, and $\D_k$ can be computed as the determinant of this matrix.

Let $\gamma: \pi \to {\rm GL}(n, R)$ be a homomorphism, where $R$ is a Noetherian UFD. One can define twisted homology of $X_\infty$ and subsequently a twisted Alexander polynomial $\D_{k,\gamma}$, well defined up to a unit in $R[t, t^{-1}]$. The classical Alexander polynomial $\D_k$ arises from the trivial homomorphism. We refer the reader to \cite{kl} for the general definition of twisted Alexander polynomial. Here we describe an equivalent way to think about $\D_{k,\gamma}$  in the case that the image of $\gamma$ is finite.

Assume that $\gamma$ is a homomorphism from $\pi$ to a finite group $\S$. Without loss of generality, we can assume that $\S$ is a group of permutation matrices in ${\rm GL}(n, {\mathbb Z})$. Let $\tilde X$ be the $n$-fold cover of $X_\infty$ determined by $\gamma$. Shapiro's Lemma implies that the ${\mathbb Z}[t, t^{-1}]$-module $H_1 \tilde X$ is isomorphic to the twisted first-homology group of $X_\infty$. It has a square presentation matrix \cite{swAGT} with determinant equal to the twisted Alexander polynomial $\D_{k,\gamma}$. 

The knot group $\pi$ acts on $\{1, \ldots, n\}$, and by restriction, so does the commutator subgroup $\pi'$. Let  ${\cal O}_\gamma$ denote the number of orbits of $\{1, \ldots, n\}$ under the action of $\pi'$. It is immediate that $H_0 \tilde X \cong {\mathbb Z}^{ {\cal O}_\gamma}$ (cf. Lemma 4.3 of \cite{fvJAMS}) and consequently its $0$th characteristic polynomial is $(t-1)^ { {\cal O}_\gamma}$ .

If $k$ is fibered, then $\D_{k,\gamma}$ is monic for any finite-image representation $\gamma$ (\cite{cha}, see also \cite{gkm}, \cite{fk}).  
A computation in section 2 of  \cite{gkm} implies that the related Wada invariant of $k$ has degree $2 g n - n$, where $g$ is the genus of $k$. By Theorem 4.1 of \cite{kl},
$$\D_{k, \gamma} = W \cdot (t-1)^{ {\cal O}_\gamma}.$$
Hence the degree of $\D_{k, \gamma}$ is 
$2 g n - n + {\cal O}_\gamma$. The following is a version of a conjecture by S. Friedl and S. Vidussi (Conjecture 4.2  \cite{fv1})  with necessary modifications for knot exteriors instead of closed $3$-manifolds. It is known for genus-1 knots \cite{fv2}.

\begin{conj}\label{folkconj} A knot $k$ is fibered iff $\D_{k,\gamma}$ is monic for every finite-image representation $\gamma$ of the knot group, and its degree is equal to $2 g n -n + {\cal O}_\gamma$.  \end{conj}

\section{Finite covers and twisted polynomials }

We recall that a group $G$ is {\it subgroup separable} if for any proper subgroup $H$ and element $g \in G\setminus H$, there exists a finite-image representation $\gamma: G \to \S$ such that $\gamma(g) \notin \gamma(H)$.

The representation $\gamma$ ``separates" the subgroup $H$ from the element $g$. The following condition can be regarded as a weak form of subgroup separability for knot groups. It requires only that the fundamental group of an incompressible Seifert surface can be separated from {\sl some} element of the commutator subgroup. 

 \begin{definition} \label{wss} The group $\pi$ of a knot $k$ is {\it weakly subgroup separable} if it admits a finite-image representation such that the image of the fundamental group of an incompressible Seifert surface for $k$ is properly contained in the image of the commutator subgroup $\pi'$. \end{definition}
 
Any incompressible Seifert surface $S$ determines an HNN decomposition $(B; U, V, \phi)$ of $\pi$ with {\it stable letter} $x$. Here the {\it base} $B$ is the fundamental group of ${\mathbb S}^3$ split along $S$, a relative cobordism between two copies of $S$, while $\phi: U=\pi_1 S \to V$ is an isomorphism. The group $\pi$ can be described as an amalgamated free product $\langle B, x \mid x^{-1} u x = \phi(u)\ (u \in U)\rangle$. The natural map $B \to \langle B, x \mid x^{-1} u x = \phi(u)\ (u \in U)\rangle$ is an embedding of $B$ into $\pi'$. (See \cite{ls} for details.)

The condition of Definition \ref{wss} can be paraphrased by saying that the {\it amalgamating subgroup} $U$ of the HNN decomposition can be separated in $\pi'$. The following result shows that the condition of Definition \ref{wss} is independent of the incompressible Seifert surface $S$, and that it is equivalent to a condition of \cite{fvJAMS}.

\begin{prop} \label{seifert} Assume that $S$ and $S'$ are incompressible Seifert surfaces for a knot $k$. If $\pi_1 S$ can be separated in $\pi'$, then it can be separated in the fundamental group of ${\mathbb S}^3$ split along $S$.  Furthermore, $\pi_1 S'$ can be separated in $\pi'$. \end{prop}

\begin{proof} Assume that $\pi_1 S$ can be separated in $\pi'$; that is, assume that there exists a finite-image representation $\gamma: \pi \to \S$ such that $\gamma(\pi_1 S)$ is a proper subgroup of $\gamma(\pi')$. 
The restriction $ \gamma|_{\pi'}$ is periodic with period (not necessarily least) equal to the order of $\gamma(x)$ in $\S$.
By Corollary 2.4(ii)  of \cite{swJPAA}, $\pi'$ has uncountably many subgroups of finite index.

Conversely, Corollary 2.4(ii) implies that if $\pi'$ has uncountably many subgroups of finite index, then for any HNN decomposition $(B; U, V, \phi)$ with stable letter $x$,  there exists a periodic finite-image representation $\rho: \pi' \to \S$ such that the sequence of subgroups $\rho(x^{-j} B x^j)$ is not constant. After conjugation by a suitable power of $x$, we can assume that some element of $\rho(B)$ is not in $\rho(x B x^{-1})$. 

We apply this first to the HNN decomposition corresponding to $S$.  Since $U=\pi_1S$ lies in $x B x^{-1}\cap B$, $\rho(\pi_1S)$ is properly contained in $\rho(B)$. 
Proposition 5.1 of \cite{swTAMS} ensures that $\rho$ extends to a representation $\bar\rho: \pi \to \bar \S$, where $\bar \S$ is a finite overgroup of $\S$, such that $\bar\rho(\pi_1S)$ is a proper subgroup of $\bar\rho(B)$. 

A similar argument, letting $(B; U, V, \phi)$ be the HNN decompostion corresponding to $S'$, shows that $\pi_1 S'$ can be separated in $B$, and hence in $\pi'$.
\end{proof}

\begin{remark} \label{Cor} Corollary 2.4(ii) of \cite{swJPAA}  and the above argument show that the group $\pi$ of a knot $k$ is weakly subgroup separable iff the infinite cyclic cover of $k$ has uncountably many finite covers.  \end{remark}

Clearly, if the  group of $k$ is weakly subgroup separable, then $k$ is nonfibered. On the other hand, if the group of every nonfibered knot $k$ is weakly subgroup separable, then Conjectures \ref{swconj} and \ref{folkconj} are true. 

A theorem of D.D. Long and G.A. Niblo \cite{ln} together with a result of D. Gabai \cite{gab2} imply that the group of any genus-1 knot is weakly subgroup separable. Friedl and Vidussi prove a more general version of this result for $3$-manifolds in \cite{fvJAMS}. See \cite{sw07} for a direct proof of our special case.

\begin{theorem}\label{main} Let $k\subset {\mathbb S}^3$ be a knot. 
The following are equivalent. 
\item{}(1) The infinite cyclic cover of $k$ has uncountably many finite covers; 

\item{}(2) the twisted Alexander polynomial $\D_{k, \gamma}$ vanishes for some finite-image representation $\gamma$; 
\item{}(3) the group of $k$ is weakly subgroup separable.
\end{theorem}

\begin{proof} We will prove that $(1) \implies (3) \implies (2) \implies (1)$. 

The implication $(1) \implies (3)$ has already been established (see Remark \ref{Cor}).

Assume that $\pi$ admits a finite-image representation $\gamma: \pi \to \S$ such that the image of an incompressible Seifert surface of $k$ is a proper subgroup of $\gamma(\pi')$. Let $n$ be the order of $\S$. If we regard $\S$ as a group of permutation matrices in ${\rm GL}(n, {\mathbb Z})$ via the right regular Cayley representation, then by Proposition \ref{seifert} and the proof of Theorem 4.2 \cite{fvJAMS}, the corresponding twisted Alexander polynomial $\D_{k, \gamma}$ vanishes. Hence $(3) \implies (2)$. 

Finally, assume that $\gamma: \pi \to \S$ is a finite-image representation such that $\D_{k,\gamma}$ vanishes. As above, we assume that $\S$ is a group of permutation matrices in ${\rm GL}(n, {\mathbb Z})$, for some $n$. Its elements permute the standard basis $\{e_1, \ldots, e_n\}$ of 
${\mathbb Z}^n$. Let $r$ be the order of $\gamma(x)$. Then the restriction $\rho= \gamma |_{\pi'}$ has period $r$; that is, $\sx^r \rho= \rho$. Consider the associated $n$-fold cover $\tilde X_\infty$ of $X_\infty$.  

The homology of $\tilde X_\infty$ is a finitely generated ${\mathbb Z}[t, t^{-1}]$-module, and its $0$th characteristic polynomial is equal to $\D_{k, \rho}=0$. We may assume that $\tilde X_\infty$ is connected; otherwise, since its homology is the direct sum of contributions from each connected component, we can replace $\tilde X_\infty$ by a connected component for which the homology has vanishing $0$th characteristic polynomial. The fundamental group of $\tilde X_\infty$ is the stabilizer of some basis vector $e_i$; that is, the set of all $a \in \pi'$ such that $\rho(a)(e_i)=e_i$.

Let $q: \hat X_\infty \to \tilde X_\infty$ be the finite $m$-fold cover such that $\pi_1 \hat X_\infty$ is the intersection $N$ of kernels 
$$N= {\rm ker}\ \rho \cap  {\rm ker}\ \sx\rho \cap \cdots
\cap  {\rm ker}\ \sx^{r-1} \rho.$$
Since $\rho$ has period $r$, the subgroup $N$ is invariant under conjugation by $x$. Hence the homology of  $H_1 \hat X_\infty$ is a finitely generated ${\mathbb Z}[t, t^{-1}]$-module.

Let $p$ a prime that does not divide $m$. Let $\tau: H_1(\tilde X_\infty; {\mathbb Z}/p)\to H_1(\hat X_\infty; {\mathbb Z}/p)$ be the transfer homomorphism induced by the chain map that takes each $i$-chain to the sum of its preimages. Since $\tau \circ q$ is multiplication by $m$, which is relatively prime to $p$, the composition $\tau \circ q$ is injective and hence so is $\tau$. Consequently,  $H_1(\hat X_\infty; {\mathbb Z}/p)$ is a submodule of $H_1(\tilde X_\infty; {\mathbb Z}/p)$, and hence the $0$th characteristic polynomial of $H_1(\hat X_\infty; {\mathbb Z}/p)$ vanishes. 

The Structure Theorem for finitely generated modules over a PID implies that $H_1(\hat X_\infty; {\mathbb Z}/p)$ contains a $({\mathbb Z}/p)[t, t^{-1}]$ summand. Consider the sequence of epimorphisms: 

$$N= \pi_1 \hat X_\infty \to H_1(\hat X_\infty; {\mathbb Z}) \to  H_1(\hat X_\infty; {\mathbb Z}/p) \to ({\mathbb Z}/p)[t, t^{-1}].$$
Regarding $({\mathbb Z}/p)[t, t^{-1}]$ as a direct sum of countably many copies of ${\mathbb Z}/p$, we see immediately that there exist uncountably many homomorphisms $h: N \to {\mathbb Z}/p$. 
Each  ${\rm ker}\ h$ is a normal in $N$ and hence normal in $\pi'$. Consider the short exact sequence
$$0 \to  {N \over {{\rm ker}\ h}} \to {\pi' \over {{\rm ker}\ h}} \to {\pi' \over N} \to 1.$$
Since $N/{{\rm ker}\ h} \cong {\mathbb Z}/p$ and $\pi'/N$ is finite, 
we see that each quotient $\pi'/{{\rm ker}\ h}$ is finite.  We can choose uncountably many $h$ such that the quotient groups 
$\pi'/{\rm ker}\ h$ are isomorphic to the same finite group $\S$.
Hence $\pi'$ has uncountably many distinct homomorphisms to $\S$, and so $(2) \implies (1)$.

\end{proof} 

%\section{Conclusion}

%We summarize the relationship between the notions discussed above. If a knot $k$ is fibered, then for every finite permutation representation $\gamma:\pi\to GL(n,\Z)$, the twisted Alexander polynomial $\D_{k,\gamma}$ is monic of degree $2 g n -n + {\cal O}_\gamma$.  In particular, $\D_{k, \gamma}\ne 0$ for all $\gamma$, which by Theorem \ref{main} is equivalent to the existence of only countably many finite covers of the infinite cyclic cover $X_\infty$.  Finally, we see from Theorem \ref{main} that this last condition implies $\gamma(\pi_1 S)$ coincides with $\gamma(\pi')$ for every finite-image representation $\gamma$ of $\pi$ and minimal-genus Seifert surface $S$.

%The weak subgroup separability condition of Definition \ref{wss} would complete the circle and give equivalence of all of the above conditions.  

%%%%%%%%%%%%%%%%%%%%%%%%%

 \bigskip

\noindent {\sl Address for both authors:} Department of Mathematics and  Statistics, ILB 325, University of South Alabama, Mobile AL  36688 USA \medskip

\noindent {\sl E-mail:} silver@jaguar1.usouthal.edu; swilliam@jaguar1.usouthal.edu

\end{document}